\documentclass[11pt]{amsart}
\usepackage{amssymb,latexsym,comment,url,amsmath}

\usepackage[T1]{fontenc}
\usepackage{rotating,graphicx}
\usepackage{xcolor}

\usepackage{diagbox} 
\usepackage{mathtools}

\usepackage{longtable}

\newcommand{\Fp}{{\mathbb{F}_p}}

\newcommand{\divv}{\operatorname{div}}

\newcommand{\Q}{{\mathbb Q}}

\newcommand{\Z}{{\mathbb Z}}

\newcommand{\Magma}{{\sf Magma}}

\newenvironment{Proof}{\par\noindent{\sc Proof:}}%
                      {\hspace*{\fill}\nobreak$\Box$\par\medskip}
                       {\hspace*{\fill}\nobreak$\Box$\par\medskip}

\newtheorem{Proposition}{Proposition}[section]
\newtheorem{Theorem}[Proposition]{Theorem}
\newtheorem{Lemma}[Proposition]{Lemma}

\newtheorem{Corollary}[Proposition]{Corollary}

\NewDocumentEnvironment{alignb}{b}{%
  \begin{align*}
  \refstepcounter{equation} #1 \tag{\theequation}
  \end{align*}
}{\ignorespacesafterend}

\allowdisplaybreaks

\theoremstyle{definition}

 \newtheorem{Example}[Proposition]{Example}
\renewcommand{\baselinestretch}{1.1}

\begin{document}

\title[Quadratic torsion orders on Jacobian varieties]%
{Quadratic torsion orders on Jacobian varieties}

\author[M. Sadek]%
{Mohammad~Sadek}
\email{mohammad.sadek@sabanciuniv.edu}
\address{Faculty of Engineering and Natural Sciences, Sabanc{\i} University, Tuzla, \.{I}stanbul, 34956 Turkey}

\author[H. Suluyer]%
{Hamide~Suluyer}
\email{hamide.kuru@atilim.edu.tr}

\address{Department of Mathematics, At{\i}l{\i}m University, 06830 G\"olba\c{s}{\i}, Ankara, Turkey}

\begin{abstract}
We establish the existence of hyperelliptic curves of genus $g\ge 2$ defined over $\mathbb{Q}$ whose Jacobians possess rational torsion points of order $N$ where $N=4g^2+2g-2$ or $4g^2+ 2g -4$. For $N = 2g^{2} + 7g + 1$, we introduce a $1$-parameter family of polynomials $f_{t}(x)$ of degree $2g+1$. For all but finitely many rational values of $t$, if the discriminant of $f_{t}(x)$ is nonzero, then the hyperelliptic curve defined by $y^{2} = f_{t}(x)$ has a rational point of order $N$ on its Jacobian. 
\end{abstract}

\maketitle

\let\thefootnote\relax\footnotetext{ \textbf{Keywords:} Hyperelliptic curves, Jacobian varieties, torsion\\

\textbf{2020 Mathematics Subject Classification:} 11G30, 14H25}

\section{Introduction}
Given an abelian variety $A$ defined over a number field $K$, the Mordell-Weil Theorem asserts that the group of rational points, $A(K)$, of $A$ is finitely generated. The subgroup $A(K)_{\operatorname{tor}}$ of torsion points of $A(K)$ is finite. In addition it is conjectured that $|A(K)_{\operatorname{tor}}|$ is uniformly bounded by a bound that depends only on the dimension $g$ of $A$ and the degree $d$ of $K$. The correctness of the conjecture was established by Merel, \cite{Merel}, when $g= 1$ and $d\ge 1$. A complete list of possible groups that can be realized as torsion subgroups when $g=1$ has been found when $d=1,2,3$, \cite{Maz78,Kam92, KM,D}. The conjecture is still wide open for any pair of integers $g\ge 2$ and $d\ge 1$.

Fixing a pair of integers $g,N\ge 2$, researchers try to construct algebraic curves over the rational field $\Q$ of genus $g$ whose Jacobian varieties possess rational torsion points of order $N$. For example, when $g=2$, infinite families of genus-$2$ curves over $\Q$ were given with torsion points on their Jacobians of order $11$ and $13$, see \cite{Bernard,Dthesis, DD, Flynn1, KS, Lep13}. Lepr\'{e}vost displayed families of genus-$2$ curves over $\Q$ with rational torsion points of order $N =15, 17,19, 21, 22, 23,24, 25, 26, 27$ and $29$ on their Jacobian varieties, see \cite{Lep2, Lep15}. In \cite{Cthesis}, the reader may find a list of integers $N$ that have appeared in the literature as orders of rational torsion points on Jacobians of algebraic curves defined over $\Q$ with genus $2,3$ or $4$.  

Flynn, \cite{Flynn1}, conjectured the existence of a constant $\kappa$, independent of $g\ge 2$, such that for every $m\leq \kappa g$, there exists a hyperelliptic curve of genus $g$ over $\Q$ with a rational $m$-torsion point on its Jacobian. Lepr\'{e}vost, \cite{LepCon}, proved this conjecture with the value $\kappa=3$. In particular, this shows the existence of abelian varieties of dimension $g$ with torsion order that is linear in $g$ for any $g\ge 2$. The latter bound has been extended to include other possible linear torsion orders in the interval $[3g,4g+1]$, see \cite{KS}. In fact, the authors show that for every integer $N$ in the interval $[3g, 4g+ 1]$,
$g \ge 3$, satisfying certain partition conditions, there exist infinite families of hyperelliptic curves
of genus $g$ whose Jacobian varieties have a rational torsion point of order $N$.

A plausible question to pose is whether it is possible to construct abelian varieties of dimension $g$ that possess rational torsion points whose order is quadratic in $g$. Fixing an even integer $g\ge 2$, Flynn gave an explicit description of $1$-parameter families of hyperelliptic curves of genus $g$ defined over $\Q$ whose Jacobian varieties contain rational torsion points of order $N$ for any $N$ in the interval $\left[g^2+2 g+1, g^2+3 g+1\right]$, see \cite{Flynn2}. For any integer $g\ge 2$, Lepr\'{e}vost displayed $1$-parameter families of hyperelliptic curves of genus $g$ with torsion points of order $2 g^2+2 g+1$ or $2 g^2+3 g+1$ on their Jacobian varieties in \cite{ Lep1}, of orders $2 g^2+4 g+1$ or $2 g(2 g+1)$ in \cite{Lep3}, and of either orders $N$, $N/2$ or $N/4$, where $N=2 g^2+5 g+5$, in \cite{Lep3}.

In this article, we extend the results of Flynn and Lepr\'{e}vost to produce  hyperelliptic curves of genus $g\ge 2$, with new rational torsion orders on their Jacobians that are quadratic in $g$. For any integer $g\ge 2$, we present hyperelliptic curves of genus $g$ 
 over $\Q$ such that their jacobian varieties contain a rational torsion point of order $N$ where $N = 4g^2+2g-2$, respectively $4g^2+ 2g -4$. Consequently, we produce the first examples in the literature of a genus-$4$ hyperelliptic curve over $\Q$ whose jacobian has a rational torsion point of order $70$, see \cite[Table 3.3]{Cthesis}. In addition, for any integer $g\ge 2$, we describe a $1$-parameter family of polynomials $f_{t}(x)$ of degree $2g+1$. For all but finitely many rational values of $t$, if the discriminant of $f_{t}(x)$ is nonzero, then the Jacobian of the hyperelliptic curve defined by $y^{2} = f_{t}(x)$ has a rational point of order $2g^{2} + 7g + 1$.
\subsection*{Acknowledgment}
The authors are indebted to the anonymous referee for the thorough reading of the manuscript and for many suggestions, comments and corrections that improved the manuscript. The authors would especially like to thank the referee for the suggestions that strengthened the statements and the proofs of Proposition \ref{prop11} and Proposition \ref{prop21}. All the calculations in this work were performed using \textbf{Magma}, \cite{Magma}. 
%This work builds upon the findings presented in the doctoral dissertation of Hamide Suluyer, \cite{Suluyer}.
  M. Sadek is supported by The Scientific and Technological Research Council of Turkey, T\"{U}B\.{I}TAK, research grant ARDEB 1001/122F312. 

\section{The construction}
\label{sec1}
Throughout this work, $K$ will be a number field. Let $f(x)\in K[x]$ be a polynomial of odd degree $2g+1$ and no repeated factors. We consider the hyperelliptic curve $C$ described by the equation $y^2=f(x)$. We let $D$ be a divisor on $C$. We recall that the {\em Riemann-Roch} space of $D$ is
the $K$-vector space $L(D)=\{\phi\in K(C): \operatorname{div}(\phi)+D\ge 0 \}\cup\{0\}.$ 

We notice that the curve $C$ has one point at infinity, $P_\infty$ since the degree of $f(x)$ is odd. There exists a natural embedding of the curve $C$ into its Jacobian $J$ that maps a point $P$ to the divisor class $[P - D]$, where $D$ is a fixed divisor of degree $1$. This map restricts to $C(K)\xhookrightarrow{} J(K)$ if $D$ is a $K$-rational divisor. In particular, one can choose $D$ to be the rational divisor $P_{\infty}$. The class of a divisor of the form $P-P_{\infty}$
defines a point in $J(K)$. We say that the divisor $P-P_{\infty}$ is a {\em torsion divisor} of order $N$ if its class in $J$ has order $N$.

We remark that $\operatorname{ord}_{P_\infty}(x)=-2$ and $\operatorname{ord}_{P_\infty}(y)=-(2g+1).$ It follows that for any $m\ge0$,  if $D_m=2(g+m+1)P_\infty$, then $L(D_m)=\langle1,x,x^2,\cdots, x^{g+m+1}, y,xy,\cdots,x^m y\rangle$; whereas if $D'_m=(2(g+m)+1)P_\infty$, then $L(D'_m)=\langle1,x,x^2,\cdots, x^{g+m}, y,xy,\cdots,x^my\rangle$.

Let $d$,  $0\le d\le g-1$, be an integer. From now on, we work with polynomials of the form  $$f(x)=A(x)^2-\lambda x^{g+1+d}(x-1)^{g-d},  \qquad A(x)\in K[x],\, \deg A(x)\le g,\quad \lambda\in K\setminus\{0\}.$$ One sees that $P_0=(0,A(0))$, ${P_0}^{\prime}=(0,-A(0))$, $P_1=(1,A(1))$, ${P_1}^{\prime}=(1,-A(1))$ are in $C(K)$. In particular, the divisors $D_i=P_i-P_\infty$ and ${D_i}^{\prime}={P_i}^{\prime}-P_\infty$, $i=0,1,$ are $K$-rational divisors on $C$.  We now consider $\phi_f\in K(C)$ such that $$\operatorname{div}(\phi_f)=(g+m){P_0}^{\prime}+(g+m+2)P_1-(2g+2m+2)P_\infty.$$

We notice that $\psi_f=\phi_f/x^{g+m}\in K(C)$ has divisor $\divv(\psi_f)=-(g+m)P_0+(g+m+2)P_1-2 P_\infty$. Therefore, one obtains the following
\begin{eqnarray}\label{eq}\left(\begin{matrix}
g+1+d&g-d\\
-(g+m)& (g+m+2)\\
\end{matrix}\right)\left(\begin{matrix}D_0\\D_1
\end{matrix}\right) = \left(\begin{matrix}\divv(y-A(x))\\ \divv(\psi_f)
\end{matrix}\right).
\end{eqnarray}

The argument above gives rise to the following result. 

\begin{Proposition}
\label{prop1}
Fix two integers $g\ge 1$ and $d$, $0\le d\le g-1$. Let $C$ be a hyperelliptic curve defined by the equation $y^2=f(x):=A(x)^2-\lambda x^{g+1+d}(x-1)^{g-d}$, where  $A(x)\in K[x]$, $\deg A(x)\le g$, $\lambda\in K\setminus\{0\}.$ We set ${P_0}^{\prime}=(0,-A(0))$, $P_1=(1,A(1)) \in C(K)$. Let $m$ be an integer such that $1\le m < d+1$. Assume that there exists a rational function $\phi_f \in K(C)$ such that
$\operatorname{div}(\phi_f)
= (g+m)P_0' + (g+m+2)P_1 - (2g+2m+2)P_\infty.$ Then there is a torsion divisor on $C$ whose order divides $2g^2+(2m+3)g+2d+m+2$. 
\end{Proposition}
\begin{Proof}
According to the argument before the proposition, the assumptions on $C$ imply the existence of a function $\psi_f=\phi_f/x^{g+m}\in K(C)$ with $\divv(\psi_f)=-(g+m)P_0+(g+m+2)P_1-2 P_\infty$. Now the result follows directly from (\ref{eq}), where $2g^2+(2m+3)g+2d+m+2$ is the determinant of the $(2\times 2)$-matrix. 
\end{Proof}

In what follows, we explain how to produce a polynomial $A(x)\in\Q[x]$ such that the order of the torsion divisor on the curve $C$ in Proposition \ref{prop1} is exactly $2g^2+(2m+3)g+2d+m+2$.

It can be seen that $\phi_f\in L((2g+2m+2)P_\infty)$. Thus, we deduce that $\phi_f=a(x)-b(x)y$ for some $a(x),b(x)\in K[x]$ such that $\deg a(x)\le g+m+1$ and $\deg b(x)\le m$. In particular, one has that the norm of $\phi_f$ in $K(C)$ is given by  \[(a(x)-b(x)y)(a(x)+b(x)y)=a(x)^2-b(x)^2y^2=h(x)x^{g+m}(x-1)^{g+m+2},\quad h(x)\in K[x].\]
Since $y^2=A(x)^2-\lambda x^{g+1+d}(x-1)^{g-d}$, it follows that
\begin{eqnarray*}
a(x)^2-b(x)^2A(x)^2&=&-\lambda b(x)^2 x^{g+1+d}(x-1)^{g-d}+h(x)x^{g+m}(x-1)^{g+m+2}\\
&=&x^{g+m}(x-1)^{g-d}\left(h(x)(x-1)^{m+2+d}-\lambda b(x)^2x^{d+1-m}\right),\quad \text{for } m<d+1.
\end{eqnarray*}

Since $P_0$ does not lie in the support of $\divv(\phi_f)$ whereas ${P_0}^{\prime}$ is in the support of $\divv(\phi_f)$, it follows that $\phi_f\left(P_0\right) \neq 0$ and $\phi_f\left({P_0 }^{\prime}\right)=0$. Moreover, $a(0)-b(0) A(0)\neq0$ whereas $a(0)+b(0) A(0)=0$. This yields that $x \nmid a(x)-b(x) A(x)$ whereas $x \mid a(x)+b(x) A(x)$. Similarly, the support of $\divv(\psi_f)$ contains $P_1$ but it does not contain ${P_1}^{\prime}$, therefore $(x-1) \mid a(x)-b(x) A(x)$ but $(x-1) \nmid a(x)+b(x) A(x)$. Thus, we may assume that 
\begin{eqnarray*}
a(x)+b(x)A(x)&=&p(x)x^{g+m},\\
a(x)-b(x)A(x)&=& q(x)(x-1)^{g-d},\\
p(x)q(x)&=&h(x) (x-1)^{m+2+d}-\lambda b(x)^2x^{d+1-m}.
\end{eqnarray*}
It follows that 
\begin{eqnarray} \label{eqnofA}
    A(x)=\frac{p(x)^2x^{g+m}-\left[h(x)(x-1)^{m+2+d}-\lambda b(x)^2x^{d+1-m}\right](x-1)^{g-d}}{2p(x)b(x)}
\end{eqnarray}
where $A(x)\in K[x]$. We set $p(x)=x-\alpha$, $\alpha\ne 0,1$, where we seek to choose $\alpha$ such that $ b(\alpha)$ and $h(\alpha)$ are non-zero, and $p(x)\mid \left( h(x)(x-1)^{m+2+d}-\lambda b(x)^2x^{d+1-m} \right)$. 

Now we will find conditions under which $b(x)$ divides the polynomial $ p(x)^2x^{g+m}-h(x)(x-1)^{g+m+2}$ of degree $m$. 
From now on, we assume that $g+m$ is even. We set $h(x)\equiv 1$. In the latter case, $p(x)^2x^{g+m}- h(x)(x-1)^{g+m+2}=(  p(x) x^{(g+m)/2}- (x-1)^{1+(g+m)/2})( p(x) x^{(g+m)/2}+ (x-1)^{1+(g+m)/2})$. Now, we set $b(x)=p(x) x^{(g+m)/2}- (x-1)^{1+(g+m)/2}$. 

 Since $p(x)=(x-\alpha) \mid\left((x-1)^{m+2+d}-\lambda b(x)^2 x^{d+1-m}\right)$, one sees that
$$
 0=(\alpha-1)^{m+2+d}-\lambda b(\alpha)^2 \alpha^{d+1-m}= (\alpha-1)^{m+2+d}-\lambda\left(-(\alpha-1)^{1+(g+m) / 2}\right)^2 \alpha^{d+1-m}.$$
In particular, one obtains that
$$ \lambda=\frac{(\alpha-1)^{m+2+d}}{(\alpha-1)^{2+g+m} \alpha^{d+1-m}} =\frac{1}{(\alpha-1)^{g-d} \alpha^{d+1-m}}.
$$

Following the discussion above, we fix an integer $m$, $1\le m< d+1$, such that $g+m$ is even. We are interested in the following $1$-parameter family of hyperelliptic curves 
\begin{eqnarray*}
C_{\alpha,d}:\; y^2&=&\left(\frac{\lambda_{\alpha,d} b_{\alpha}(x)x^{d+1-m}(x-1)^{g-d}+ (x-1)^{1+(g+m)/2}+(x-\alpha) x^{(g+m)/2}}{2(x-\alpha)}\right)^2 -\lambda_{\alpha,d} x^{g+1+d}(x-1)^{g-d},\\
A_{\alpha,d}(x)&=&\left(\frac{\lambda_{\alpha,d} b_{\alpha}(x)x^{d+1-m}(x-1)^{g-d}+ (x-1)^{1+(g+m)/2}+(x-\alpha) x^{(g+m)/2}}{2(x-\alpha)}\right),\\
b_{\alpha}(x)&=&  (x-\alpha) x^{(g+m)/2}-(x-1)^{1+(g+m)/2},\quad \alpha\ne 0,1,\\
%%\lambda&=&\frac{(\alpha-1)^{1+(g+m)/2}}{b_{\alpha}(\alpha) \alpha^{g+1-m}}.
\displaystyle \lambda_{\alpha,d}&=&\frac{ 1}{ \alpha^{d+1-m}(\alpha-1)^{g-d}}.
\end{eqnarray*}
It is clear that $\deg(b_{\alpha}(x))=(g+m)/2$ if $\alpha \neq 1+(g+m)/2$,
whereas $\deg(b_{\alpha}(x))=-1+(g+m)/2$ if $\alpha = 1+(g+m)/2$.
The latter observation imposes certain restrictions on the possible values of $\alpha$ and $d$. We recall that $\deg(b_{\alpha}(x))$ is at most $m$. If $\deg(b_{\alpha}(x))=(g+m)/2$, then $g\le m$ which contradicts the fact that $m<d+1\le g$. If $\deg(b_{\alpha}(x))=-1+(g+m)/2$, then $m \geq g - 2$ and the assumption that $g+m$ is even implies that $m=g-2$, hence $\alpha=g$ and $d\in\{g-1,g-2\}$.

Besides the above choice of $b(x)$, we will also consider the case where $b(x)$ is chosen to be a linear polynomial, hence $m=1$.
In this case, the parameters $\lambda$ and $\alpha$ are selected so that $A(x)$ becomes a polynomial. This will be treated in detail in \S~5. 

 %One sees that $\deg (b(x))=(g+m)/2$ if $\alpha\neq (g+m)/2+1$, whereas $\deg (b(x))=(g+m)/2-1$ otherwise. Since $\deg (b(x))$ is at most $m$ and $\deg A(x)=g$, the possible values for $m$ are either $g$ and then $\alpha\neq{g+1}$; or $g-2$ where in this case $\alpha=g$. The fact that $m<d+1\le g$ implies that $m=g-2$, hence $\alpha=g$ and $b (x)=(x-g)x^{g-1}-(x-1)^g$. 

\section{The genus of the hyperelliptic curves}
In this section, we show that the curve $C_{\alpha,d}$ defined in the previous section where $\alpha=g$ and $d$ is either $g-1$ or $g-2$ is indeed a hyperelliptic curve of genus $g$ over $\Q$.

We recall that the {\it height} $h(f)$ of a polynomial $f\in\Z[x]$ is the height of its
coefficients, namely, if $f(x)=a_dx^d+\cdots+a_0$, then $h(f)=\log\max\{1,|a_i|\}$, see \cite[B.7]{Hindry}. If $\alpha$ is a non-zero algebraic number of degree $d$ with minimal polynomial $f(x)=a_d\prod_{i=1}^d(x-\alpha_i)$ over $\Z$ with $a_d>0$, we define the {\it Mahler measure} of $\alpha$, $M(\alpha)$, by $M(\alpha)=a_d\prod_{i=1}^d\max\{1,|\alpha_i|\}$. The height of $\alpha$, $h(\alpha)$, is defined by $h(\alpha)=\log M(\alpha)/d$.
\begin{Proposition}
\label{prop11}
 Fix an integer $g\ge 2$. We set $\lambda_{g,g-1}=1/(g^2(g-1))$. The curve $C_{g,g-1}$ defined over $\Q$ by 
\begin{eqnarray*}
C_{g,g-1}&:&y^2= f_{g-1}(x)= A_{g,g-1}(x)^2 -4\lambda_{g,g-1}\,x^{2g}(x-1),
\end{eqnarray*}
where 
\begin{eqnarray*}
A_{g,g-1}(x)&=& \frac{(x-g) x^{g-1}+(x-1)^g+\lambda_{g,g-1} \,x^2(x-1)\left((x-g) x^{g-1}-(x-1)^g\right)}{(x-g)},
\end{eqnarray*}
is a hyperelliptic curve of genus $g$.
\end{Proposition}
\begin{Proof}
 We need to prove that $f_{g-1}(x)$ is a square-free polynomial. We assume on the contrary that $f_{g-1}(x_0)=f_{g-1}'(x_0)=0$ for some root $x_0$ of $f_{g-1}(x)$ in the algebraic closure $\overline{\Q}$ of $\Q$. 

We may write \begin{eqnarray}\label{expr}A_{g,g-1}(x)=x^{g-1}p_{g-1}(x)-(x-1)^gq_{g-1}(x),\end{eqnarray} where $$p_{g-1}(x)=1+\lambda_{g,g-1}x^2(x-1),\qquad q_{g-1}(x)=\lambda_{g,g-1}(x^2+(g-1)x+g^2 -g).$$ In particular, $x_0\ne0,1$. One now sees that

\begin{eqnarray*}
f_{g-1}^{\prime}(x)&=&2 A_{g,g-1}(x) A_{g,g-1}^{\prime}(x)-8g\lambda_{g,g-1} x^{2 g-1}(x-1)-4\lambda_{g,g-1} x^{2 g} \\
&=& 2 A_{g,g-1}(x) A_{g,g-1}^{\prime}(x)-4\lambda_{g,g-1} x^{2 g-1}((2g+1)x-2g).
\end{eqnarray*}
The vanishing of both $f_{g-1}(x)$ and $f^{\prime}_{g-1}(x)$ at $x_0$ yields that
\begin{eqnarray*}
2 x_0\left(x_0-1\right) A_{g,g-1}^{\prime}\left(x_0\right)=\left((2 g+1) x_0-2 g\right) A_{g,g-1}\left(x_0\right).
\end{eqnarray*}
Substituting (\ref{expr}) in the latter identity gives rise to the following 
\begin{eqnarray}\label{exp2}
x_0^{g-1} G_{g-1}\left(x_0\right)=(x_0-1)^g H_{g-1}\left(x_0\right), 
\end{eqnarray}
  where $H_{g-1}(x), G_{g-1}(x) \in \mathbb{Q}[x]$ are defined by 
{\footnotesize\begin{eqnarray*}
 G_{g-1}(x)&=&(-3 x+2)p_{g-1}(x)+2 x(x-1) p_{g-1}^{\prime}(x)=3 \lambda_{g,g-1} x^4-5\lambda_{g,g-1} x^3+2 \lambda_{g,g-1} x^2-3x+2\\&=&q_{g-1}(x)(x-g)(3x-2),\\
H_{g-1}(x)&=&(-x+2 g)q_{g-1}(x)+2 x(x-1) q_{g-1}^{\prime}(x)=3 \lambda_{g,g-1}x^3+(3 g \lambda_{g,g-1}-5 \lambda_{g,g-1}) x^2+\left(g^2 \lambda_{g,g-1}-3 g \lambda_{g,g-1}+2 \lambda_{g,g-1}\right) x  +2.
\end{eqnarray*}}

The resultant of $G_{g-1}(x)$ and $H_{g-1}(x)$ as polynomials in $x$ is given by
$$\operatorname{Res}_x(G_{g-1}(x),H_{g-1}(x))=
72\lambda_{g,g-1}^7  (3g + 1)^2  (g - 1)^3 g^5  (9g^2 - 10 g + 2).
$$
Therefore, $\operatorname{Res}_x(G_{g-1}(x),H_{g-1}(x))$ is non-zero when $g\ge 2$. In other words, $G_{g-1}(x)$ and $H_{g-1}(x)$ have no common roots in $\overline{\Q}$. In particular, using (\ref{exp2}), one sees that
$H_{g-1}\left(x_0\right) G_{g-1}\left(x_0\right) \neq 0$ for any $x_0 \neq 0,1$, and for any $g \geq 2$. More precisely, 
 one has $x_0^{g-1} G_{g-1}\left(x_0\right)/H_{g-1}\left(x_0\right)=\left(x_0-1\right)^g$.

We define the polynomial $S_{g-1}(x)$ as follows %and $T_{g-1}(x)$ as follows
{\footnotesize\begin{eqnarray*}
 S_{g-1}(x)&:=& (H_{g-1}(x) p_{g-1}(x)-G_{g-1}(x) q_{g-1}(x))^2-4 \lambda_{g,g-1} x^2(x-1) H_{g-1}(x)^2 \\
 &=&-4\lambda_{g,g-1}(x-g)q_{g-1}(x)^2 M_{g-1}(x),%\\
%T_{g-1}(x)&:=&2\left( H_{g-1}(x)p_{g-1}(x)-G_{g-1}(x) q_{g-1}(x)\right)\Big((g-1) H_{g-1}(x)(x-1) p_{g-1}(x)+x(x-1) H_{g-1}(x) p_{g-1}^{\prime}(x)\\&-&g xq_{g-1}(x) G_{g-1}(x)-x(x-1) G_{g-1}(x) q_{g-1}^{\prime}(x)\Big)
%-4\lambda_{g,g-1}  x^2(x-1) H_{g-1}(x)^2\left(2 g(x-1)+ x\right)\\
%&=&-4\lambda_{g,g-1}(x-g)((2g+1)x-2 g )q_{g-1}(x)^2 M_{g-1}(x).\\
\end{eqnarray*}
}
where $M_{g-1}(x)$ is the polynomial defined by
$$M_{g-1}(x)=9 x^4 + 4(2g-5 ) x^3 + 2(3g^2 - 8 g + 
    7) x^2 + 3 ( g-1)^3 x + g(g-1)^3 .$$
Evaluating $S_{g-1}(x)$ %and $T_{g-1}(x)$
at $x_0$ gives the equality
$S_{g-1}(x_0)=H_{g-1}^2(x_0) f_{g-1}(x_0)/x_0^{2 g-2}$. %and $T_{g-1}(x_0)=H_{g-1}^2(x_0)(x_0-1) f_{g-1}^{\prime}(x_0)/x_0^{2 g-3}$.
In particular, a common root $x_0$ of $f_{g-1}(x)$ and $f_{g-1}^{\prime}(x)$ is a root of $S_{g-1}(x)$, %and $T_{g-1}(x)$, 
i.e., $x_0$ is either $g$, a root of $q_{g-1}(x)$, or a root of $M_{g-1}(x)$. However, one recalls that $G_{g-1}(x)=q_{g-1}(x)(x-g)(3x-2)$. Therefore, $x_0$ is neither $g$ nor a root of $q_{g-1}(x)$ since otherwise $x_0$ would be a root of $G_{g-1}(x)$ which contradicts (\ref{exp2}) and the fact that $G_{g-1}(x)$ and $H_{g-1}(x)$ do not have common roots. Now we are left with showing that none of the roots of $M_{g-1}(x)$ is a common root of $f_{g-1}(x)$ and $f_{g-1}^{\prime}(x)$. 
We assume on the contrary that $x_0$ is a multiple root of $f_{g-1}(x)$ and that $M_{g-1}(x_0)=0$. According to (\ref{exp2}), one has 
\begin{eqnarray*}
\alpha^g=x_0 H^*_{g-1}(x_0)/G^*_{g-1}(x_0),
\end{eqnarray*}
where
\begin{eqnarray*}
\alpha=x_0/(x_0-1), \quad G_{g-1}^*(x)=\lambda_{g,g-1}^{-1}G_{g-1}(x)\in\Z[g][x],\quad H_{g-1}^*(x)=\lambda_{g,g-1}^{-1}H_{g-1}(x)\in\Z[g][x]. 
\end{eqnarray*}
We first claim that $\alpha$ and $\alpha^{-1}$ are not conjugates, in particular, $\alpha$ is not a root of unity. We notice that $\alpha$ is a root of the polynomial 
{\footnotesize\begin{eqnarray*}
M_1(x):=(x-1)^4M_{g-1}(x/x-1)= g^4x^4&+&(-4 g^4 +3g^3+3g^2+g+1)x^3+(6g^4-9g^3-3g^2+5g+5)x^2\\&-&(4g+3)(g-1)^3x+g(g-1)^3,
\end{eqnarray*}}
whereas $\alpha^{-1}$ is a root of the polynomial 
    $M_2(x):=x^4M_1(1/x)$. The resultant of $M_1(x)$ and $M_2(x)$ is given by $$\operatorname{Res}_x(M_1(x),M_2(x))=9 (-2 + g)^2 g^2 (-1 + 2 g + 7 g^2 - 12 g^3 + g^4)^2 (1 + 8 g - 
   24 g^3 + 16 g^4),$$ where the degree-4 factors are irreducible in $\Q[g]$. Therefore, for any integer $g>2$, $\alpha$ and $\alpha^{-1}$ are not conjugates. In particular, $\alpha$ is not a root of unity.
   
For $g\ge 3$, a result of Smyth together with the fact that $\alpha$ and $\alpha^{-1}$ are not conjugates imply that $M(\alpha)>\theta=1.32471\cdots$, where $\theta$ is the real root of the polynomial $x^3-x-1$, see \cite{Smyth}. It follows that 
\begin{eqnarray}
\label{eq:h}
h(\alpha^g)=gh(\alpha)\ge g\log(1.32471)/4\approx 0.0702984\cdot g.
\end{eqnarray}
In addition, one sees that 
\begin{eqnarray*}
    h\left(x_0 H^*_{g-1}(x_0)/G^*_{g-1}(x_0)\right) &\le& h(x_0)+h(H^*_{g-1}(x_0))+h(G^*_{g-1}(x_0))\\
    &\le&8h(x_0)+h(H^*_{g-1})+h(G^*_{g-1})+\min\{\log 4,4\log 2\}+\min\{\log 5,5\log 2\},
\end{eqnarray*} where the second inequality follows from \cite[Part B, Proposition B.7.1]{Hindry}. Landau's inequality, \cite[Chapter 3, Proposition 2.7]{Lang} implies that 
\begin{eqnarray*}
h(x_0)\le \frac{1}{8}\log\left(81+16(2g-5)^2+4(3g^2-8g+7)^2+9(g-1)^6+g^2(g-1)^6\right)\le \frac{1}{8}\log (2801g^8).
\end{eqnarray*}
Moreover, one has $h(G^*_{g-1})=\log(3g^2(g-1))$ and $h(H^*_{g-1})=\log(2g^2(g-1))$. Therefore, 
$$h\left(x_0 H^*_{g-1}(x_0)/G^*_{g-1}(x_0)\right)\le \log (16806 g^{14})+\log20\le 13+14 \log(g). $$
Combining the latter inequality with the inequality (\ref{eq:h}), one has $13+14\log (g)\ge 0.0702984\cdot g.$ It follows that if $x_0$ is a multiple root of $f_{g-1}(x)$ and a root of $M_{g-1}(x)$, then one must have that $g\le 1662.$ One may check using \Magma, \cite{Magma}, that $C_{g,g-1}$ is indeed a hyperelliptic curve for all $g$, $2\le g\le 1662$, or equivalently that the discriminant of $f_{g-1}(x)$ is non-zero for these values of $g$. Thus, the proof is concluded for $C_{g,g-1}$.
\end{Proof}

We now prove that $C_{g,g-2}$ is a hyperelliptic curve of genus $g$.

\begin{Proposition}
\label{prop21}
 Fix an integer $g\ge 2$. We set $\lambda_{g,g-2}=1/(g(g-1)^2)$. The curve $C_{g,g-2}$ defined over $\Q$ by 
\begin{eqnarray*}
C_{g,g-2}&:& y^2= f_{g-2}(x)=A_{g,g-2}(x)^2 -4\lambda_{g,g-2}\,x^{2g-1}(x-1)^2,
\end{eqnarray*}
where 
\begin{eqnarray*}
A_{g,g-2}(x)&=& \frac{(x-g) x^{g-1}+(x-1)^g+\lambda_{g,g-2}\, x(x-1)^2\left((x-g) x^{g-1}-(x-1)^g\right)}{(x-g)},
\end{eqnarray*}
is a hyperelliptic curve of genus $g$.
\end{Proposition}
\begin{Proof}
 We now need to prove that $f_{g-2}(x)$ is a square-free polynomial. We assume on the contrary that $f_{g-2}(x_0)=f_{g-2}'(x_0)=0$ for some root $x_0$ of $f_{g-2}(x)$ in the algebraic closure $\overline{\Q}$ of $\Q$. 

We may write \begin{eqnarray}\label{expr-2}A_{g,g-2}(x)=x^{g-1}p_{g-2}(x)-(x-1)^gq_{g-2}(x),\end{eqnarray} where $$p_{g-2}(x)=1+\lambda_{g,g-2}x(x-1)^2,\qquad q_{g-2}(x)=\lambda_{g,g-2}(x^2+(g-2)x+(g-1)^2).$$ In particular, $x_0\ne0,1$. One now sees that

\begin{eqnarray*}
f_{g-2}^{\prime}(x)&=&2 A_{g,g-2}(x) A_{g,g-2}^{\prime}(x)-4\lambda_{g,g-2} x^{2 g-2}(x-1)((2g-1)(x-1)+2x).
\end{eqnarray*}
The evaluation of both $f_{g-2}(x)$ and $f^{\prime}_{g-2}(x)$ at $x_0$ yields that
\begin{eqnarray*}
2 x_0\left(x_0-1\right) A_{g,g-2}^{\prime}\left(x_0\right)=\left((2 g-1) x_0+2 x_0\right) A_{g,g-2}\left(x_0\right).
\end{eqnarray*}
Substituting (\ref{expr-2}) in the latter identity gives rise to the following 
\begin{eqnarray}\label{exp3}
x_0^{g-1} G_{g-2}\left(x_0\right)=(x_0-1)^g H_{g-2}\left(x_0\right), 
\end{eqnarray}
  where $H_{g-2}(x), G_{g-2}(x) \in \mathbb{Q}[x]$ are defined by 
{\footnotesize\begin{eqnarray*}
 G_{g-2}(x)&=&(-3 x+1)p_{g-2}(x)+2 x(x-1) p_{g-2}^{\prime}(x)=\lambda_{g,g-2}(
3x^4 - 7x^3 + 5x^2 - (3g^3 - 6g^2 + 3g + 1)x)  +1 \\&=&q_{g-2}(x)(x-g)(3x-1),\\
H_{g-2}(x)&=&(-x+2 g-1)q_{g-2}(x)+2 x(x-1) q_{g-2}^{\prime}(x)= \lambda_{g,g-2} (3x^3+(3g - 7)x^2   + (g^2 - 5g + 5)x)+2-1/g.
\end{eqnarray*}}

The resultant of $G_{g-2}(x)$ and $H_{g-2}(x)$ as polynomials in $x$ is given by
$$\operatorname{Res}_x(G_{g-2}(x),H_{g-2}(x))=
72 \,\lambda_{g,g-2}^7\, (3g -4)^2  
g^3 (g-1)^5  (9g^2 - 8 g + 1).
$$
Therefore, $\operatorname{Res}_x(G_{g-2}(x),H_{g-2}(x))$ is non-zero when $g\ge 2$. In other words, $G_{g-2}(x)$ and $H_{g-2}(x)$ have no common roots in $\overline{\Q}$. In particular, using (\ref{exp3}), one sees that
$H_{g-2}\left(x_0\right) G_{g-2}\left(x_0\right) \neq 0$ for any $x_0 \neq 0,1$, and for any $g \geq 2$. More precisely, 
 one has $x_0^{g-1} G_{g-2}\left(x_0\right)/H_{g-2}\left(x_0\right)=\left(x_0-1\right)^g$.
We define the polynomial $S_{g-2}(x)$ as follows %and $T_{g-2}(x)$ as follows
{\footnotesize\begin{eqnarray*}
 S_{g-2}(x)&:=& (H_{g-2}(x) p_{g-2}(x)-G_{g-2}(x) q_{g-2}(x))^2-4 \lambda_{g,g-2} x(x-1)^2 H_{g-2}(x)^2 \\
 &=&-4\lambda_{g,g-2}(x-g)q_{g-2}(x)^2 M_{g-2}(x),\\
%T_{g-2}(x)&:=&2\left( H_{g-2}(x)p_{g-2}(x)-G_{g-2}(x) q_{g-2}(x)\right)\Big((g-1) H_{g-2}(x)(x-1) p_{g-2}(x)+x(x-1) H_{g-2}(x) p_{g-2}^{\prime}(x)\\&-&g xq_{g-2}(x) G_{g-2}(x)-x(x-1) G_{g-2}(x) q_{g-2}^{\prime}(x)\Big)
%-4\lambda_{g,g-2}  x(x-1)^2 H_{g-2}(x)^2\left((2 g-1)(x-1)+ 2x\right)\\
%&=&-4\lambda_{g,g-2}(x-g)((2g+1)x-2 g+1 )q_{g-2}(x)^2 M_{g-2}(x)
\end{eqnarray*}}
where $M_{g-2}(x)$ is defined as follows
$$M_{g-2}(x)=
9x^4 + 8(g - 3)x^3 + 2(3g^2 - 10g + 11)x^2 + (3g^3 - 12g^2 + 16g - 8)x + (g - 1)^4. $$
Evaluating $S_{g-2}(x)$ %and $T_{g-2}(x)$
at $x_0$ gives
$S_{g-2}(x_0)=H_{g-2}^2(x_0) f_{g-2}(x_0) / x_0^{2g-4}$. %and
%$T_{g-2}(x_0)=H_{g-2}^2(x_0)(x_0-1) f_{g-2}'(x_0) / x_0^{2g-3}$.
In particular, any common root $x_0$ of $f_{g-2}(x)$ and $f_{g-2}'(x)$ is also a root of $S_{g-2}(x)$. % and $T_{g-2}(x)$.
Thus, $x_0$ is either $g$, or a root of $q_{g-2}(x)$, or a root of $M_{g-2}(x)$.
Recall that $G_{g-2}(x)=q_{g-2}(x)(x-g)(3x-1)$.
Hence, $x_0$ cannot be equal to $g$ nor a root of $q_{g-2}(x)$, since otherwise $x_0$ would be a root of $G_{g-2}(x)$, contradicting $(\ref{exp3})$ and the fact that $G_{g-2}(x)$ and $H_{g-2}(x)$ have no common roots.
Therefore, it remains to show that no root of $M_{g-2}(x)$ is a common root of $f_{g-2}(x)$ and $f_{g-2}'(x)$.

We assume on the contrary that $x_0$ is a multiple root of $f_{g-2}(x)$ and that $M_{g-2}(x_0)=0$. According to (\ref{exp3}), one has 
\begin{eqnarray*}
\beta^g=x_0 H^*_{g-2}(x_0)/G^*_{g-2}(x_0),
\end{eqnarray*}
where
\begin{eqnarray*}
\beta=x_0/(x_0-1), \quad G_{g-2}^*(x)=\lambda_{g,g-2}^{-1}G_{g-2}(x)\in\Z[g][x],\quad H_{g-2}^*(x)=\lambda_{g,g-2}^{-1}H_{g-2}(x)\in\Z[g][x]. 
\end{eqnarray*}
The algebraic integers $\beta$ and $\beta^{-1}$ are not conjugates, in particular, $\beta$ is not a root of unity. This can be checked as follows. One has $\beta$ is a root of the polynomial 
{\footnotesize\begin{eqnarray*}
N_1(x):=(x-1)^4M_{g-2}(x/x-1)= g^3(g-1)x^4&-&g^3(4g-7)x^3+(6g^4-15g^3+6g^2+4g+4)x^2\\&+&(-4g^4+13g^3-12g^2+4)x+(g-1)^4,
\end{eqnarray*}}
whereas $\beta^{-1}$ is a root of the polynomial 
    $N_2(x):=x^4N_1(1/x)$. The resultant of $N_1(x)$ and $N_2(x)$ is given by $$\operatorname{Res}_x(N_1(x),N_2(x))=9 (-1 + g)^2 (1 + g)^2 (-3 + 16 g - 23 g^2 + 8 g^3 + g^4)^2 (1 + 
   24 g^2 - 40 g^3 + 16 g^4),$$ where the degree-4 factors are irreducible in $\Q[g]$. Therefore, for any integer $g\ge2$, $\beta$ and $\beta^{-1}$ are not conjugates.

   We now use the same argument we used for the curve $C_{g,g-1}$. We have 
   \begin{eqnarray*}
   0.0702984\cdot g\le h(\beta^g)&\le& 8h(x_0)+h(H^*_{g-2})+h(G^*_{g-2})+\log 20\\
   &\le & \log (5186g^8)+\log(2g^3)+\log(13g^3)+\log20=14\log(g)+15.
   \end{eqnarray*} 
It follows that $g\le 1695.$ One may check using \Magma, \cite{Magma}, that $C_{g,g-2}$ is indeed a hyperelliptic curve for all $g$, $2\le g\le 1695$.
 Thus, the statement of the proposition holds for $C_{g,g-2}.$ 
\end{Proof}

\section{The order of the torsion divisor}
In this section, we discuss the order of the torsion subgroup of the Jacobian of the curves $C_{g,g-1}$ and $C_{g,g-2}$ defined in the previous sections.
From now on we write $a_{d}(x)$ and $b_{d}(x)$ for $a(x)$ and $b(x)$ that were introduced in \S \ref{sec1} to help keep track of the parameter $d$, where $d$ is either $g-1$ or $g-2$.
\begin{Theorem}
\label{thm1}
Fix an integer $g\ge 2$. We set $\lambda_{g,g-1}=1/(g^2(g-1))$ and $\lambda_{g,g-2}=1/(g(g-1)^2)$. We consider hyperelliptic curves of genus $g$ defined over $\Q$ by the following equations
\begin{eqnarray*}
C_{g,g-1}&:&y^2= f_{g-1}(x)= A_{g,g-1}(x)^2 -4\lambda_{g,g-1}x^{2g}(x-1),\\
C_{g,g-2}&:& y^2= f_{g-2}(x)=A_{g,g-2}(x)^2 -4\lambda_{g,g-2}x^{2g-1}(x-1)^2,\\
\end{eqnarray*}
where 
\begin{eqnarray*}
A_{g,g-1}(x)&=& \frac{(x-g) x^{g-1}+(x-1)^g+\lambda_{g,g-1} x^2(x-1)\left((x-g) x^{g-1}-(x-1)^g\right)}{(x-g)},\\
A_{g,g-2}(x)&=& \frac{(x-g) x^{g-1}+(x-1)^g+\lambda_{g,g-2} x(x-1)^2\left((x-g) x^{g-1}-(x-1)^g\right)}{(x-g)}.
\end{eqnarray*}
The torsion divisor $D_0=P_0-P_{\infty}$ on the curve $C_{g,g-1}$, respectively the torsion divisor $D_1=P_1-P_{\infty}$ on the curve $C_{g,g-2}$, has order  $4g^2+2g-2$, respectively $4g^2+2g-4$.
\end{Theorem}
\begin{Proof}
That the curves $C_{g,g-1}$ and $C_{g,g-2}$ are of genus $g$ over $\Q$ follows from Proposition \ref{prop11} and Proposition \ref{prop21}. 
Let $a_{g-1}$ and $b_{g-1}$ be given as follows:
$$
\begin{aligned}
&b_{g-1}(x)= (x-g) x^{g-1}-(x-1)^{g}\\
&a_{g-1}(x)=(x-g)x^{2g-2}-b_{g-1}(x)A_{g,g-1}(x)
\end{aligned}$$
We recall the existence of the following rational functions on the curve $C_{g,g-1}$

$$\begin{aligned}
& {\phi_{f_{g-1}}}=a_{g-1}(x)-b_{g-1}(x) y, \\
& \psi_{f_{g-1}}=\frac{{\phi_{f_{g-1}}}}{x^{g+m}}=\frac{a_{g-1}(x)-b_{g-1}(x) y}{x^{g+m}}=\frac{a_{g-1}(x)-b_{g-1}(x) y}{x^{2g-2}}, \\
& {\theta_{f_{g-1}}}=y-A_{g,g-1}(x),
\end{aligned}
$$
where the norm of ${\phi_{f_{g-1}}}$ is given by $a_{g-1}^2(x)-b_{g-1}^2(x) y^2=x^{g+m}(x-1)^{g+m+2}=x^{2g-2}(x-1)^{2g}$.
According to Proposition \ref{prop1}, the order of the class of the divisor $D_0=(0,A_{g,g-1}(0))-P_{\infty}$ divides $l=4 g^2+2 g-2$. It follows that the principal divisor $l D_0$ is the divisor of the rational function $L_{g-1}(x, y)$ where
\begin{eqnarray}\label{eqn1}
 L_{g-1}(x, y)=\frac{{{\theta_{f_{g-1}}}} ^{g+m+2}}{{\psi_{f_{g-1}}}^{g-d}}
 =\frac{{\theta_{f_{g-1}}}^{2 g} \cdot x^{2 g-2}}{{\phi_{f_{g-1}}}}=\frac{(y-A_{g,g-1}(x))^{2g}x^{2g-2}}{a_{g-1}(x)-b_{g-1}(x)y}. 
\end{eqnarray}

\newcommand\numberthis{\addtocounter{equation}{1}\tag{\theequation}}

This implies that
\begin{align*}
 L_{g-1}(x, y)
 =& \frac{(y-A_{g,g-1}(x))^{2 g} \cdot(a_{g-1}(x)+b_{g-1}(x) y) \cdot x^{2 g-2}}{a_{g-1}(x)^2-b_{g-1}(x)^2 y^2} \\
  =&\frac{(y-A_{g,g-1}(x))^{2 g} \cdot(a_{g-1}(x)+b_{g-1}(x) y) \cdot x^{2 g-2}}{(x-1)^{2 g} \cdot x^{2 g-2}}   \\
  =&\frac{(y-A_{g,g-1}(x))^{2 g} \cdot(a_{g-1}(x)+b_{g-1}(x) y)}{(x-1)^{2 g}} \numberthis\label{eqn2}
  \end{align*}
   % =&\frac{(y-A_{g,g-1}(x))^{2 g} \cdot(A_{g,g-1}(x)+b_{g-1}(x) y)}{(x-1)^{2 g}} \cdot \frac{(y+A_{g,g-1}(x))^{2g}}{(y+A_{g,g-1}(x))^{2g}} \\
 % =&\frac{\left(y^2-A_{g,g-1}^2(x)\right)^{2g} \cdot(A_{g,g-1}(x)+b_{g-1}(x) y)}{(x-1)^{2 g} (y+A_{g,g-1}(x))^{2g}} 
% \\ =&\frac{\left(-4\lambda_{g,g-1} x^{2 g}(x-1)\right)^{2 g} \cdot(A_{g,g-1}(x)+b_{g-1}(x) y)}{(x-1)^{2 g} (y+A_{g,g-1}(x))^{2g}} \\
 % =&\frac{(-4\lambda_{g,g-1})^{2 g} \cdot x^{4g^2} \cdot(A_{g,g-1}(x)+b_{g-1}(x) y) }{(y+A_{g,g-1}(x))^{2g}}\numberthis\label{eqn3}

We recall that ${P_1}^{\prime}=(1,-A_{g,g-1}(1)), {P_0}^{\prime}=(0,-A_{g,g-1}(0)) \in C_{g,g-1}(\Q)$ do not appear in the support of the divisor $D_0$. From (\ref{eqn1}) and (\ref{eqn2}), we can compute $L_{g-1}\left({P_1}^{\prime}\right)$ and $ L_{g-1}\left({P_0}^{\prime}\right)$ as follows
\begin{eqnarray*}
L_{g-1}\left({P_1}^{\prime}\right) & =&\frac{(-2  A_{g,g-1}(1))^{2 g} \cdot 1^{2 g-1}}{a_{g-1}(1)+b_{g-1}(1)  A_{g,g-1}(1)}  =\frac{2^{2 g}}{(1-g)}, \\
L_{g-1}\left({P_0}^{\prime}\right) & =&\frac{(-2 A_{g,g-1}(0))^{2 g} \cdot(a_{g-1}(0)-b_{g-1}(0) A_{g,g-1}(0))}{(-1)^{2 g}}=\frac{2^{2g} \cdot}{(-g)^{2 g} \cdot g\cdot(-1)^{2g-1}}=\frac{-2^{2 g}}{g^{2 g+1}}. 
\end{eqnarray*}
% L_g\left(P_1\right)&=\frac{(-4\lambda_{g,g-1})^{2 g} \cdot(a_g(1)+b_g(1)  A_{g,g-1}_{g-1}(1))}{{(2 A_{g,g-1}_{g-1}(1))}^{2g}}=\frac{ 2^{2g}}{(g-1)^{2g} \cdot g^{4g}}(1-g)  =\frac{2^{2g}}{(1-g)^{2g-1} \cdot g^{4g}}.}

%In addition by equation \ref{test} we can obtain that 
%\begin{eqnarray*}

%\end{eqnarray*} 
%$$
%\begin{aligned}
%& L\left({P_1}^{\prime}\right)=a N^m\left({P_1}^{\prime}\right)=\frac{-1}{g^{2 g+1} } \\
%& L\left(P_1\right)=a N^m\left(P_1\right)=\frac{-1}{(g-1)^{2 g-1} \cdot g^{4g}}
%\end{aligned}

 We therefore obtain the following identity
$$L_{g-1}\left({P_1}^{\prime}\right)(g-1)=L_{g-1}\left({P_0}^{\prime}\right)  g^{2g+1}.$$
%%Given that $2g+1$ is relatively prime to $l=4g^2+2g-2=2(2g-1)(g+1)$, this implies that the order of the class of $D_0$ cannot be a proper divisor of $l$, hence it must be $l$ itself. 
If the order of $D_0$ is $t$, then $l=t\cdot s$ for some $s$. This implies the existence of a rational function $N_{g-1}$ on the curve $C_{g,g-1}$ such that $tD_0$ is the divisor of $N_{g-1}$. Given that $ lD_0=st D_0$ is the divisor of $L_{g-1}$, we obtain that $L_{g-1}=uN_{g-1}^s$, for some $u\in \mathbb{Q}\setminus\{0\}$. In particular, we get
$$N^s_{g-1}\left({P_1}^{\prime}\right)(g-1)=N^s_{g-1}\left({P_0}^{\prime}\right)  g^{2g+1}.$$
It can be checked that $g^{2g+1}/(g-1)$ can not be an $s$-th power for any $s>1$ that is relatively prime to $2g+1$. Given that $s$ is a divisor of 
  $l=4g^2+2g-2=2(2g-1)(g+1)$ and that $2g+1$ is relatively prime to $l$, the order of the class of $D_0$ cannot be a proper divisor of $l$, and hence must be $l$ itself.

As for the curve $C_{g,g-2}$, we set 
$b_{g-2}(x)= (x-g) x^{g-1}-(x-1)^{g},\ a_{g-2}(x)=(x-g)x^{2g-2}-b_{g-2}(x)A_{g,g-2}(x)
$, ${\phi_{f_{g-2}}}=a_{g-2}(x)-b_{g-2}(x) y$, $\psi_{f_{g-2}}=\phi_{f_{g-2}}/x^{g+m}=\phi_{f_{g-2}}/x^{2g-2}$, and
$ \theta_{f_{g-2}}=y-A_{g,g-2}(x)$. We consider the class of the divisor $D_1:=P_1-P_{\infty}$, where $P_1=(1,A_{g,g-2}(1))$. According to Proposition \ref{prop1}, there is a rational function $L_{g-2}(x, y)$ defined on $C_{g,g-2}$ such that the principal divisor $l' D_1$ is the divisor of $L_{g-2}$, where $l'=4 g^2+2 g-4$. In fact, the function $L_{g-2}$ is defined as follows

\begin{align*}
L_{g-2}(x, y) & =\theta_{f_{g-2}}^{2 g-2} \cdot \psi_{f_{g-2}}^{2 g-1}=\frac{\theta_{f_{g-2}}^{2 g-2} \cdot \phi_{f_{g-2}}^{2 g-1}}{x^{(2 g-1) \cdot(2 g-2)}} \\
& =\frac{(y-A_{g,g-2}(x))^{2 g-2} \cdot(a_{g-2}(x)-b_{g-2}(x) y)^{2 g-1}}{x^{(2 g-1) \cdot(2 g-2)}} \numberthis\label{eqn6}\\
& =\frac{(y-A_{g,g-2}(x))^{2 g-2} \cdot(a_{g-2}(x)-b_{g-2}(x) y)^{2 g-1} \cdot(a_{g-2}(x)+b_{g-2}(x) y)^{2 g-1}}{x^{(2 g-1)(2 g-2)} \cdot(a_{g-2}(x)+b_{g-2}(x) y)^{2 g-1}} \\
& =\frac{(y-A_{g,g-2}(x))^{2 g-2} \cdot\left(a_{g-2}^2(x)-b_{g-2}^2(x) y^2\right)^{2 g-1}}{x^{(2 g-1)(2 g-2)} \cdot(a_{g-2}(x)+b_{g-2}(x) y)^{2 g-1}} \\
& =\frac{(y-A_{g,g-2}(x))^{2 g-2} \cdot x^{(2 g-2) \cdot(2 g-1)} \cdot(x-1)^{(2 g-4)(2 g-1)}}{x^{(2 g-2)(2 g-1)} \cdot(a_{g-2}(x)+b_{g-2}(x) y)^{2 g-1}} \\
& =\frac{(y-A_{g,g-2}(x))^{2g-2} \cdot(x-1)^{(2 g-4) \cdot(2 g-1)}}{(a_{g-2}(x)+b_{g-2}(x) y)^{2 g-1}} \numberthis\label{eqn7} \\
\end{align*}
For the points ${P_1}^{\prime}=(1,-A_{g,g-2}(1)), {P_0}^{\prime}=(0,-A_{g,g-2}(0)) \in C_{g,g-2}(\Q)$, we use (\ref{eqn6}) and (\ref{eqn7}) to compute  $ L_{g-2}\left({P_1}^{\prime}\right)$ and $L_{g-2}\left(P_0^{\prime}\right)$ respectively as follows
\begin{align*}
 L_{g-2}\left({P_1}^{\prime}\right) & =\frac{(-2 A_{g,g-2}(1))^{2 g-2} \cdot(1-g)^{2 g-1}}{1^{(2 g-1) \cdot(2 g-2)}} =2^{2g-2}(1-g)^{2 g-1},\\
 L_{g-2}\left({P_0}^{\prime}\right)&=\frac{(-2 A_{g,g-2}(0))^{2 g-2} \cdot(-1)^{(2 g-4) \cdot(2 g-1)}}{(a(0)-b(0) A_{g,g-2}(0))^{2 g-1}}=\frac{2^{2 g-2} \cdot\left(\frac{1}{g}\right)^{2 g-2}}{\left(\frac{-1}{g}\right)^{2 g-1}}  =-2^{2 g-2} \cdot g .
\end{align*}
%& =\frac{1}{(1-g)^{2 g-2} \cdot g^{4 g{-3}}}\frac{\operatorname{Lg}\left(P_0^{\prime}\right)}{\lg \left(P_1^{\prime}\right)}=\frac{-2^{2 g-2} \cdot g}{2^{2 g^{-2}} \cdot(1-g)^{2 g-1}} \Rightarrow g N^m\left(P_1^{\prime}\right)=(1-g)^{2 g-1} N^m\left(P_0^{\prime}\right)

%& =\frac{(y-A(x))^{2 g-2} \cdot(y+A(x))^{2 g-2} \cdot(a(x) \cdot b(x) y)^{2 g-1}}{x^{(2 g-1) \cdot(2 g-2)}(y+A(x))^{2 g-2}} \\
%& =\frac{\left(-\lambda \cdot x^{2 g-1} \cdot(x-1)^2\right)^{2 g-2} \cdot(a(x)-b(a) y)^{2 g-1}}{x^{(2 g-1) \cdot(2 g-2)}(y+A(x))^{2 g-2}} \\
%& =\frac{(-\lambda)^{2 g-2} \cdot(x-1)^{2 \cdot(2 g-2)} \cdot(a(x)-b(x) y)^{2 g-1}}{(y+A(x))^{2 g-2}}

%$$
%\begin{aligned}
%L\left(P_0\right) & =\frac{\left(\frac{1}{(g-1)^2 \cdot g}\right)^{2 g-2} \cdot(-1)^{2 \cdot(2 g-2)} \cdot q(0)^{2 g^{-1}}}{\left(\frac{(-1)^g}{(1-g)}\right)^{2 g-2}} \\

%& =\frac{(1-g)^{2 g{-2}}}{(g-1)^{4 g{-4}} \cdot g^{2 g-2} \cdot g^{2 g-1}} \\
%& =\frac{1}{(1-g)^{2 g-2} \cdot g^{4 g{-3}}}
%\end{aligned}
%$$

 We notice that $g L_{g-2}\left({P_1}^{\prime}\right)=(g-1)^{2 g-1} L_{g-2}\left({P_0}^{\prime}\right)$. It can be easily verified that $(g-1)^{2g-1}/g$ is not an $s$-th power for any $s>1$ that is relatively prime to $2g-1$. These facts together with the observation that $2g-1$ is relatively prime to $4g^2+2g-4=2(2g-1)(g+1)-2$ imply that the order of the class of the divisor of $D_1$ is exactly $4g^2+2g-4$. 
\end{Proof}
In the following table, we produce hyperelliptic curves of genus $g$, $2\le g\le 5$, whose Jacobians possess rational torsion points with order determined by Theorem \ref{thm1}. The curves $C_{g,g-1}$ and $C_{g,g-2}$ appear in the table as $C_g$ and $C'_g$, respectively.

{\centering
{\footnotesize\begin{table}[h]
\label{Tab}
\begin{tabular}{|c|c|c|c|}
\hline
$g$ & \textbf{Curve} & \textbf{Torsion divisor} & \textbf{Order} \\
\hline
$2$ & 
$ C_2:y^2  = -16 x^5 + 17 x^4 - 14 x^3 
    +53 x^2 - 28 x + 4$
& $D_0 = (0 : 2 : 1) - (1 : 0 : 0)$ & $18$ \\
\hline
$2$ & 
$C_2':y^2 = -8 x^5 + 17 x^4 - 16 x^3+
    18 x^2 - 8 x + 1$
 & $D_1 = (1 : 2 : 1) - (1 : 0 : 0)$ & $16$ \\
\hline
$3$ & $C_3: y^2 =  -72  x^7 + 81  x^6 - 186  x^5+
    1057  x^4 - 1028  x^3 + 628  x^2 - 192  x + 36$  & $D_0 = (0 : 6 : 1) - (1 : 0 : 0)$ & $40$ \\
\hline
$3$ & $C_3': y^2 =  -48x^7 + 105x^6 - 180x^5 +
    550x^4 - 508x^3 + 297x^2 - 88x + 16$  & $D_1 = (1 : 12 : 1) - (1 : 0 : 0)$ & $38$ \\
\hline
$4$ & $\begin{aligned}C_4:y^2 &=  -192  x^9 + 228  x^8 - 984  x^7 +
    7456  x^6 - 10544  x^5 + 11245  x^4\\ -& 7458  x^3 + 3489  x^2 - 1080  x + 144 \end{aligned}$  & $D_0 = (0 : 12 : 1) - (1 : 0 : 0)$ & $70$ \\
\hline
$4$ & $\begin{aligned}C_4':y^2 &=  -144x^9 + 324x^8 - 912x^7 +
    4660x^6 - 6424x^5 + 6669x^4 \\&- 4348x^3 + 2002x^2 - 612x + 81 \end{aligned}$ & $D_1 = (1 : 36 : 1) - (1 : 0 : 0)$ & $68$ \\
\hline
$5$ & $\begin{aligned}C_5: y^2 &= -400 x^{11} + 500 x^{10} - 3400 x^9+ 32200 x^8 - 59720 x^7 + 90685 x^6\\&- 92770 x^5 + 71241 x^4 - 41352 x^3+
    16456 x^2 - 3840 x + 400\end{aligned}$  & $D_0 = (0 : 20 : 1) - (1 : 0 : 0)$ & $108$ \\
\hline
$5$ & $\begin{aligned}
    C_5':y^2 &=  -320x^{11} + 740x^{10} - 3120x^9
    + 22300x^8- 40720x^7 + 60645x^6 \\ &- 61300x^5 + 46586x^4 - 26804x^3 +
    10601x^2 - 2464x + 256
\end{aligned}$ & $D_1 = (1 : 80 : 1) - (1 : 0 : 0)$ & $106$ \\
\hline
\end{tabular}
\caption{Curves of low genus and different torsion orders}
\end{table}}}
An abelian variety defined over a field $K$ is called $K$-simple if it is not isogenous over $K$ to a product of abelian varieties
of lower dimensions. 
\begin{Corollary}
The Jacobian varieties of the curves in Table 1 are $\Q$-simple varieties. 
\end{Corollary}
\begin{Proof}
We refer the reader to \cite{water} for the information introduced in this paragraph. Let $p$ be a good prime of a hyperelliptic curve $C$. We define the Zeta function of $C$ by $\displaystyle Z_C (t) = \operatorname{exp}\left(\sum_{n\ge 1}|C(\mathbb{F}_{p^n})|\frac{t^n}{n}\right).$ In fact, $Z_C(t)$ can be written in the
form $L_C(t)/(1 - t)(1 - pt)$
where $L_C (t) \in \Z[t]$ is of degree $2g$. Moreover, 
$L_C (t) = t^{
2g}P_C (1/t)$, where $P_C (t)$ is the characteristic polynomial of the Frobenius endomorphism of the Jacobian of $C$. If the Jacobian of $C$
is $\Fp$-simple, then it can be shown that $P_C(x) = h(x)^e$ where $h(x)\in  \Z[x]$ is irreducible over $\Z$ and $e \ge 1$. 

In what follows, we choose a prime $p$ of good reduction for each of the curves in Table 1. For $C_2$, we choose the prime $p=31$; whereas the prime $p=23$ is a prime of good reduction for $C_i$, $i=3,4,5$ and $C_j'$, $j=2,3,4,5$. Moreover, one may use \Magma\, to find $L_{C_i}(t)$ and $L_{C_i'}$ over $\Fp$. In addition, it can be verified that these polynomials $L_C(t)$ are irreducible. In other words, the Jacobian varieties of the curves $C_i$ and $C_i'$ are $\Fp$-simple. If the Jacobian of a curve $C$ is not $\Q$-simple, then the Jacobian of $C$ as an abelian variety over $\Fp$ is not $\Fp$-simple, where $p$ is a prime of good reduction for $C$. This concludes the proof. 
\end{Proof}
As pointed out by the referee, in view of \cite[Table 3.1]{Cthesis}, there is a hyperelliptic curve of genus $2$ whose Jacobian variety contains a rational torsion point of order $70$. This implies the existence of a $4$-dimensional non-simple abelian variety with a rational torsion point of order $70$. The Jacobian of the curve $C_4$ is a $\Q$-simple abelian variety with a torsion divisor of order $70$.

 \section{An infinite family of hyperelliptic curves}
In Proposition \ref{prop1}, for a fixed integer $g\ge 2$, one may choose $m=1$ and $1\le d\le g-1$. Under the assumption of existence of the function $\phi$ with $\operatorname{div}(\phi)=(g+1)P_0^\prime+(g+3)P_1-(2g+4)P_\infty$, it follows that there exists a hyperelliptic curve of genus $g$ with a divisor whose class represents a rational point in the Jacobian of order dividing $2 g^2+5 g+2 d+3$.
 We fix two integers $g$ and $d$ with $g\ge 3$ and $0\le d\le g-1$. Following the construction in \S \ref{sec1}, we may give explicit $1$-parameter families of such curves for which the order of the divisor is exactly $2g^2+5g+2d+3$ when $d$ is chosen to be $g-1$ and $h(x)=1$.

\begin{Theorem}
\label{thm2}
Fix an odd integer $g\ge 3$. We set
\begin{eqnarray*}
\alpha=\beta-\frac{(\beta-1)^{\frac{g+3}{2}}}{\beta^{\frac{g+1}{2}}},\qquad
\beta\in \mathbb{Q}\setminus\{0,1\},
\qquad
\lambda=\frac{(\alpha-1)^{g+2}}{(\alpha-\beta)^{2}\alpha^{g-1}}\in\mathbb{Q}.
\end{eqnarray*}
We also set
$$
A(x)=\frac{x^{g+1}(x-\alpha)^{2}-(x-1)\left((x-1)^{g+2}-\lambda(x-\beta)^{2}x^{g-1}\right)}{2(x-\alpha)(x-\beta)},
\qquad
y^2=f(x)=A(x)^2-\lambda x^{2g}(x-1).
$$

Then $A(x)$ is a polynomial of degree $g$ in $\mathbb{Q}[x]$. In addition, for all but finitely many rational values of $\beta$, if the discriminant of $f(x)$ is nonzero, then the equation $
y^2=f(x)$
defines a hyperelliptic curve $C_\beta$ of genus $g$ over $\mathbb{Q}$, where the divisor
$D_0=(0,A(0))-P_\infty$
on $C_\beta$ is a torsion divisor whose order is exactly $2g^2+7g+1$.
\end{Theorem}

% \begin{Theorem}
%\label{thm2}
%Fix an odd integer $g\ge 3$. We set 
%\begin{eqnarray*}
%\alpha=\beta-\frac{(\beta-1)^{\frac{g+3}{2}}}{\beta^ \frac{g+1}{2}},\, \beta\in \mathbb{Q}-\{ 0,1\},\quad \lambda=\frac{(\alpha-1)^{g+2}}{(\alpha-\beta)^{2} \alpha^{g-1}} \in \Q.
%\end{eqnarray*}
%We also set \[A(x)=\frac{x^{g+1}(x-\alpha)^{2}-(x-1)\left((x-1)^{g+2}-\lambda(x-\beta)^{2} x^{g-1}\right)}{2(x-\alpha)(x-\beta)},\qquad y^2=f(x)=A(x)^2-\lambda x^{2g}(x-1).\]

 %Then $A(x)$ is a polynomial of degree $g$ in $\Q[x]$. In addition, if the discriminant of $f(x)$ is nonzero, then the equation 
%$
%y^{2}=f(x)$ defines a hyperelliptic curve $C_{\beta}$ of genus $g$ over $\Q$, where
 %the divisor $D_0=(0,A(0))-P_{\infty}$ on $C_{\beta}$ is a torsion divisor whose order is $2g^2+7g+1$.
%\end{Theorem}
\begin{Proof}
 %For a fixed value of $g$, the discriminant of $f(x)$ is a polynomial in $\beta$. Therefore, there are at most finitely many values of $\beta$ such that the discriminant of $f(x)$ vanishes. Thus, the $1$-parameter equation defines a family of hyperelliptic curves of genus $g$ over $\Q(\beta)$.
Let $B(x)=x^{g+1}(x-\alpha)^{2}-(x-1)\left((x-1)^{g+2}-\lambda(x-\beta)^{2} x^{g-1}\right)$. One sees that $B(\alpha)=-(\alpha-1)((\alpha-1)^{g+2}
 -\lambda(\alpha-\beta)^2\alpha^{g-1})$. Substituting with the expression for $\lambda$, one gets $B(\alpha)=0$. Similarly, $B(\beta)=0$, hence $(x-\alpha)(x-\beta)$ divides $B(x)$ and $A(x)\in\Q[x]$. In addition, $\deg A(x)=g$, therefore, $\deg f(x)=2g+1$.
 
In view of Proposition \ref{prop1}, the order of the class of the divisor $D_0$ divides $l=2g^2+7g+1$. Let $L(x,y)$ be a rational function such that $l D_0$ is the divisor of $L(x,y)$.
We set $\theta(x,y)=y-A(x)$, $\bar{\theta}(x,y)=y+A(x)$,  $\phi(x,y)=a(x)-b(x) y$, $\bar{\phi}(x,y)=a(x)+b(x) y$ and $\psi(x,y)=\phi(x,y)/x^{g+1}$,  where $a(x)$ and $b(x)$ are as in \S \ref{sec1} satisfying the following identities 
$$
    a(x)-b(x) A(x)=(x-1) q(x),\qquad
    a(x)+b(x) A(x)=x^{g+1} p(x),
$$
where $p(x) q(x)=(x-1)^{g+2}-\lambda b^2(x) x^{g-1}$, $p(x)=x-\alpha$ and $b(x)=x-\beta$. Then we obtain $a(x)=x^{g+1}(x-\alpha)-(x-\beta)A(x)$. The function $L(x,y)$ may be defined as follows 
$$
L(x, y) =\frac{\theta(x,y)^{g+3}}{\psi(x,y)} =\frac{\theta(x,y)^{g+3}}{\phi(x,y)} x^{g+1} %=\frac{x^{g+1} \cdot \theta(x,y)^{g+2}(\theta(x,y) \bar{\varphi})}{\varphi \cdot \bar{\varphi}}
$$ %where
%$\bar{\varphi}=a(x)+b(x) y$ and  $\varphi\cdot\bar{\varphi}=a^{2}(x)-b^{2}(x) y^{2} =x^{g+1}(x-1)^{g+3} 
%$.
%If follows that
%$$
%L(x, y)=\frac{(y-A(x))^{g+2}\left[q(x)(y-A(x))-\lambda x^{2 g}(x-\beta)\right]}{(x-1)^{g+2}},\qquad q(x)=\frac{(x-1)^{g+2}-\lambda(x-\beta)^{2} x^{g-1}}{(x-\alpha)}. \label{1}
%$$
If the order of the class of $D_0$ is not $l$, then we may assume that $l=m n$, where the order of the class of $D_{0}$ is $n$. In particular, there exists a rational function $N(x, y)$ such that the principal divisor $nD_0$ is the divisor of $N(x,y)$. It follows that there exists $u\in\Q\setminus\{0\}$ such that 
 $u N^{m}(x, y)=L(x, y)$. Since $P_{1}^{\prime}=(1,-A(1))$ and $P_{0}^{\prime}=(0, -A(0))$ are not in the support of the divisor $D_{0}$, one obtains that
$
L(P_{1}^{\prime})/L(P_{0}^{\prime})=N^{m}(P_{1}^{\prime})/N^{m}(P_{0}^{\prime}) $. 
Evaluating the function $L(x,y)$ at $P_1'$, one gets
\begin{eqnarray*}L\left(P_{1}^{\prime}\right)=\cfrac{(-2 A(1))^{g+3}}{a(1)+b(1)A(1)}=\cfrac{(1-\alpha)^{g+2}}{(\beta-1)^{g+3}}.
\end{eqnarray*}

Similarly, one sees that
\begin{align*}
L(x,y)
&=\frac{(y-A(x))^{g+3}\,x^{g+1}(a(x)+b(x)y)}
{a(x)^2-b(x)^2y^2}=\frac{(y-A(x))^{g+3}\,x^{g+1}(a(x)+b(x)y)}
{x^{g+1}(x-1)^{g+3}}\\
&=\frac{(y-A(x))^{g+3}(a(x)+b(x)y)}{(x-1)^{g+3}}.
\end{align*}
Therefore, 
\begin{eqnarray*}
L(P_0')&=(2A(0))^{g+3}\,(a(0)-b(0)A(0))
= \frac{-1\cdot(-1)^{(g+2)(g+3)}}{\alpha^{g+3}\beta^{g+3}(-\alpha)}= -\frac{1}{\alpha^{g+4}\beta^{g+3}}.
\end{eqnarray*}
Then we obtain that
\[
\frac{L(P_1')}{L(P_0')}
=
\frac{(\alpha-1)^{g+2}\,\alpha^{g+4}\,\beta^{g+3}}
{(\beta-1)^{g+3}}.
\]
Since $g$ is odd, we assume that $g=2k-1$ for some integer $k\ge 2$.
We know that
$$
\alpha
= \beta - \frac{(\beta-1)^{k+1}}{\beta^k}= \frac{\beta^{k+1} - (\beta-1)^{k+1}}{\beta^k}.
$$

Hence,
$$
\alpha - 1
= \frac{\beta^{k+1} - \beta^k - (\beta-1)^{k+1}}{\beta^k}= (\beta-1)\left(\frac{\beta^k - (\beta-1)^k}{\beta^k}\right).
$$

Therefore,
$$
\frac{L(P_1')}{L(P_0')}
=
\frac{
(\beta-1)^{g+2}
(\beta^k-(\beta-1)^k)^{g+2}
(\beta^{k+1}-(\beta-1)^{k+1})^{g+4}
\beta^{g+3}
}{
\beta^{k(g+2)}\,
\beta^{k(g+4)}\,
(\beta-1)^{g+3}
}.$$

Clearing denominators, we obtain
$$
L(P_0')\,
(\beta^k-(\beta-1)^k)^{g+2}\,
(\beta^{k+1}-(\beta-1)^{k+1})^{g+4}
=
L(P_1')\,\beta^{4k^2+2k-2}(\beta-1).
$$
$$N^{m}(P_{0}^{\prime})\,
(\beta^k-(\beta-1)^k)^{g+2}\,
(\beta^{k+1}-(\beta-1)^{k+1})^{g+4}
=
N^{m}(P_{1}^{\prime})\,\beta^{4k^2+2k-2}(\beta-1).$$
The square-free polynomials
$\beta^k-(\beta-1)^k,\quad
\beta^{k+1}-(\beta-1)^{k+1},\quad
\beta,\quad
\beta-1
$ are pairwise relatively prime.
We note that
\[
\gcd(g+2,\,2g^2+7g+1)=
\begin{cases}
5, & \text{if } g \equiv 3 \pmod{5},\\
1, & \text{otherwise}.
\end{cases}
\]

Similarly,
\[
\gcd(g+4,\,2g^2+7g+1)=
\begin{cases}
5, & \text{if } g \equiv 1 \pmod{5},\\
1, & \text{otherwise}.
\end{cases}
\]
It follows that at least one of the integers $g+2$ or $g+4$ is relatively prime to $2g^2+7g+1$. Let $q\ge 2$ be a prime divisor of $m$ such that either $\gcd(q,g+2)=1$ or $\gcd(q,g+4)=1$. One may see that $
\frac{
(\beta^k-(\beta-1)^k)^{g+2}
(\beta^{k+1}-(\beta-1)^{k+1})^{g+4}}
{\beta^{g(g+3)}(\beta-1)}
$
can not be a $q$-power except possibly for finitely many rational values of $\beta$. This can be justified as follows. One considers the following equation $$y^{q}=(\beta^k-(\beta-1)^k)^{a}\,
(\beta^{k+1}-(\beta-1)^{k+1})^{b}\beta^{c}(\beta-1)^{q-1}$$ where $g+2\equiv a$, $g+4\equiv b$ and $-g(g+3)\equiv c$ modulo $q$, $0\le a,b,c\le q-1$. Since either $\gcd(q,g+2)=1$ or $\gcd(q,g+4)=1$, it follows that $\max(a,b)\ge1$. For $q>2$, the latter equation describes a superelliptic curve whose genus can be computed by the Riemann–Hurwitz
formula as in \cite[Proposition 3.7.3]{Henning} yielding that the genus is at least $2$. For \(q=2\) and $g\ge 5$, the curve is hyperelliptic with genus at least $2$.
In both cases, Faltings' theorem implies the existence of finitely many rational points on this curve. 
The remaining case is $q=2$ and $g=3$, where the latter curve becomes $y^2=6x^4 - 15x^3 + 14x^2 - 6x + 1$. This curve can be described by the Weierstrass equation $y^2 = x^3 - x^2 + x$ whose Mordell-Weil rank can be checked using \Magma \, \cite{Magma} to be 0.
This argument together with the fact that $l$ is relatively prime with either $g+2$ or $g+4$ imply that the order of $D_0$ is exactly $l$ for all but finitely many rational values of $\beta$.
\end{Proof}

\begin{Example}
 Consider the genus-$3$ hyperelliptic curve described by $y^2 =  -190512 x^7 + 727801 x^6 -
    1181596 x^5 + 1054252 x^4 - 527008 x^3 + 166448 x^2 - 30912 x + 3136$. The class of the divisor  $D_0=(0:56:1)-(1:0:0)$ is of order $40$ in the Jacobian of the hyperelliptic curve. This curve corresponds to $\beta=2$ and $g=3$ in Theorem \ref{thm2}. We remark that this curve is not isomorphic to the genus-$3$ curve $C_3$ given in Table 1.  

Consider the genus-$5$ hyperelliptic curve defined by
$y^2=43740000 x^{11}-21753279 x^{10}-231167700 x^9+1056603628 x^8-2750786720 x^7+4528772176 x^6-4990905280 x^5+3807828800 x^4-2015296000 x^3+717280000 x^2-156800000 x+16000000$.
%$y^2  = -1334139660000 x^{11} +7810111072849 x^{10} - 21237895715004 x^9 + 35509499052172 x^8 -40961959673568 x^7 + 34980786907216 x^6 - 22600090770240 x^5 +11170251259200 x^4 - 4105259712000 x^3 + 1045133280000 x^2 - 162518400000 x+ 11664000000$. 
The class of the divisor $D_0=(0:4000:1)-(1:0:0)$ gives rise to a rational point of order $86$ in the Jacobian of the curve. This curve corresponds to $\beta=2/3$ and $g=5$ in Theorem \ref{thm2}.
\end{Example}

Although we stated Theorem \ref{thm2} for curves of odd genus, the same technique extends to curves of even genus as follows.

\begin{Theorem}
\label{cor}
In Theorem \ref{thm2}, we set $\beta=u(t):=\frac{(t^2+1)^2}{4t^2}$
where $t\neq0,\pm1$ is a rational number. 
 For all but finitely many rational values of $t$,
if the curve $C_{u(t)}$ is a hyperelliptic curve of genus $g$ defined over $\Q$, then it possesses a torsion divisor of order $2g^2+7g+1$.
\end{Theorem}

\begin{Proof}
    In Theorem \ref{thm2}, the expression $\beta=u(t)$ forces both $\beta$ and $\beta-1$ to be rational squares. This allows $\alpha$ to be rational for any choice of an integer $g\ge 2$. Now the proof follows by adjusting the proof of Theorem \ref{thm2}.
\end{Proof}

\begin{Example} Consider the genus-$2$ hyperelliptic curve described by
$y^2 = -299054816676000 x^5 +
    937313042871529 x^4 - 1165161421194050 x^3 + 677279473485625 x^2 -
    132825168000000 x + 8294400000000$. The class of $D_0=(0:2880000:1)-(1:0:0)$ is of order $23$. This curve corresponds to $t=2$ and $\beta=25/16$ in Corollary \ref{cor}.\\
    
    \noindent Consider the genus-$4$ hyperelliptic curve described by
$ y^2 =
    -441076451313968208343861667771372100000 x^9\\
    +2231009503403670702562982043605865222649 x^8 -
    4959972109544667027708192318400142478050 x^7 +
    6329054704630532302814017899191191335625 x^6 -
    5260199601304122072610634289700416000000 x^5 +
    3123070596609213073858989244272000000000 x^4 -
    1315926242281486797139217238210000000000 x^3 +
    396345305950692328102018752000000000000 x^2 -
    76786692290915614316668800000000000000 x +\\
    6733822227188480528040000000000000000$. The class of $D_0=(0:2594960929800000000:1)-(1:0:0)$ is of order $61$. This curve corresponds to $t=2$ and $\beta=25/16$.
    
\end{Example}

\section*{Conflict of Interest statement} On behalf of all authors, the corresponding author states that there is no conflict of interest.
\section*{Data availability statement} The authors declare that the data supporting the findings of this study are available within the paper.

\end{document}